\documentclass[12pt, a4]{amsart}

\usepackage[a4paper, left=28mm, right=28mm, top=28mm, bottom=34mm]{geometry}
\usepackage[all]{xy}
\usepackage{amsmath}
\usepackage{amssymb}
\usepackage{amsthm}
\usepackage{mathrsfs}
\theoremstyle{definition}
\newtheorem{thm}{Theorem}[section]

\theoremstyle{definition}

\numberwithin{equation}{section}

\title[Algorithms of linear determinantal representations of plane curves]{On algorithms to obtain linear determinantal representations of smooth plane curves of higher degree}

\author{Yasuhiro Ishitsuka}
\address{Department of Mathematics, Faculty of Science, Kyoto University, Kyoto 606-8502, Japan}
\email{yasu-ishi@math.kyoto-u.ac.jp}

\author{Tetsushi Ito}
\address{Department of Mathematics, Faculty of Science, Kyoto University, Kyoto 606-8502, Japan}
\email{tetsushi@math.kyoto-u.ac.jp}

\author{Tatsuya Ohshita}

\address{Graduate School of Science and Engineering, Ehime University, Matsuyama, 790-8577, Japan}
\email{ohshita.tatsuya.nz@ehime-u.ac.jp}

\date{\today}
\subjclass[2010]{Primary 11D41; Secondary 14H50, 14K15}
\keywords{linear determinantal representations, plane quartics, the Klein quartic, the Fermat quartic}

\def\Q{{\mathbb Q}}

\def\Z{{\mathbb Z}}

\def\P{{\mathbb P}}

\DeclareMathOperator{\Homsheaf}{\mathscr{H}\!\textit{om}}

\def\Br{\mathop{\mathrm{Br}}\nolimits}

\def\Aut{\mathop{\mathrm{Aut}}\nolimits}

\def\Jac{\mathop{\mathrm{Jac}}\nolimits}
\def\Pic{\mathop{\mathrm{Pic}}\nolimits}

\def\id{\mathop{\mathrm{id}}\nolimits}

\newcommand{\transp}[1]{{}^{t}\!{#1}}

\newtheorem{Alg}[thm]{Algorithm}			

\begin{document}

\maketitle

\section{Introduction}
Let $C \subset \P^2$ be a smooth plane curve of degree $d \ge 1$ defined over a field $k$.
It is defined by a single homogeneous polynomial 
\[
    F = F(X,Y,Z) \in k[X, Y, Z]
\]
of degree $d$.
\textit{A linear determinantal representation of $C$ over $k$} is a square matrix $M = (l_{i,j})$ of size $d$ 
whose entries $l_{i,j} = l_{i,j}(X,Y,Z)$ are $k$-linear forms in the variables $X, Y, Z$ satisfying
\[
	C = \left\{ 
		[X:Y:Z] \in \P^2 \mid \mathrm{det}(M) = 0
	\right\} \subset \P^2.
\]
Equivalently, the determinant of $M$ coincides with $F$ as a polynomial in $X, Y, Z$ up to a non-zero constant.

Two linear determinantal representations $M, M'$ are said to be \textit{equivalent} if
there exist two invertible matrices $A, B \in \mathrm{GL}_d(k)$ such that
\(
	M' = AMB.
\)
If $M$ is symmetric, we say that $M$ is \textit{a symmetric linear determinantal representation}. Two symmetric linear determinantal representations are said to be equivalent if
they are equivalent as linear determinantal representations.

Studying linear determinantal representations is a classical topic in algebraic geometry {\cite{Bea00, Dol12}}.
Recently it also appears in the study of derived categories \cite{Gal14}, 
semi-definite programming \cite{PSV11}. 
There are also studies from arithmetic viewpoints \cite{FN14, Ish2, II1, II2}.

In this paper, we continue and extend the first author's results in \cite{Ish1}, where an algorithm computing linear determinantal representations of smooth plane cubics was given. 
Here we shall give two algorithms to compute linear determinantal representations explicitly for smooth plane curves of any degree. 
As particular examples, we give all linear determinantal representations, up to equivalence, of 
the Klein quartic and the Fermat quartic over the field $\Q$ of rational numbers.

\section{Linear determinantal representations and line bundles}\label{Param}

Let $k$ be a field, $d \ge 1$ a positive integer, and $C \subset \P^2$ a smooth plane curve of degree $d$.
Its genus $g$ is equal to $(d-1)(d-2)/2$.
The canonical bundle $\omega_C$ of $C$ is
isomorphic to the Serre twist $\mathcal{O}_C(d-3)$.
A line bundle $\mathcal{L}$ on $C$ is said to be \textit{non-effective} if
it has no non-zero global sections, i.e.\ $H^0(C, \mathcal{L}) = 0$.
A \textit{theta characteristic} on $C$ is a line bundle $\mathcal{L}$ on $C$ satisfying $\mathcal{L}^{\otimes 2} \cong \omega_C$.
The degree of a theta characteristic is equal to $g-1$.

The following classical theorem is a foundation of our algorithms.

\begin{thm}\label{Bij}
	Let $C \subset \P^2$ be a smooth plane curve of degree $d$ defined over $k$.
	\begin{enumerate}
    	\item There is a natural bijection between the following two sets:
    	\begin{itemize}
    		\item the set of equivalence classes of linear determinantal representations of $C$ over $k$, and
    		\item the set of isomorphism classes of non-effective line bundles on $C$ of degree $g - 1 = d(d-3)/2$.
    	\end{itemize}
		In this bijection, the operation of transpose $M \mapsto \transp{M}$ corresponds to
		the operation 
		\[
		    \mathcal{L} \mapsto \mathcal{L}^{\vee} := \Homsheaf_{\mathcal{O}_C}(\mathcal{L}, \mathcal{O}_C(d-3)),
		\]
		where $\Homsheaf_{\mathcal{O}_C}$ is
		the sheaf Hom between $\mathcal{O}_C$-modules.
    	\item This bijection induces a bijection between the following two subsets:
    	\begin{itemize}
    		\item the set of equivalence classes of symmetric linear determinantal representations of $C$ over $k$, and
    		\item the set of isomorphism classes of non-effective theta characteristics on $C$.
    	\end{itemize}
	\end{enumerate}
\end{thm}

For the proof of this theorem and related results, see \cite{Bea00, Dol12, Bea77}. For plane cubics, see also \cite{Ish1}. 

We briefly explain the construction of the bijection in Theorem \ref{Bij} (1).
Let us take a non-effective line bundle $\mathcal{L}$ of degree $g-1$ on the curve $C$. 
Let $R = k[X,Y,Z]$ be the homogeneous coordinate ring of $\P^2$. 
It is well-known that the graded $R$-module
\begin{align*}
    N &:= \Gamma_* (C, \mathcal{L}) = \bigoplus_{n \in \Z} H^0(C, \mathcal{L}(n))
\end{align*}
has a minimal free resolution of the form
\begin{equation*}\label{Eq: LDR}
	\xymatrix@C=10pt{
		0 \ar[r] & R(-2) \otimes_k W_1 \ar[r]^{\widetilde{M}} & 
		R(-1) \otimes_k W_0 \ar[r] & N \ar[r] & 0.
	}
\end{equation*}
Here $W_0, W_1$ are $d$-dimensional $k$-vector spaces defined by
\begin{align*}
    W_0 &:= H^0(C, \mathcal{L}(1)), \\
    W_1 &:= \ker \big(
    m \colon H^0(C, \mathcal{L}(1)) \otimes_k H^0(C, \mathcal{O}_C(1)) \\
    &\qquad \qquad \to H^0(C, \mathcal{L}(2))
    \big).
\end{align*} 
After fixing bases of these vector spaces, 
the $R$-homomorphism $\widetilde{M}$ can be written as a square matrix $M$ of 
size $d$ whose entries are $k$-linear forms in the variables $X, Y, Z$. 
This matrix $M$ gives a linear determinantal representation corresponding to $\mathcal{L}$. It can be shown that the equivalence class of $M$ is uniquely determined by the isomorphism class of the line bundle $\mathcal{L}$, and every equivalence class of linear determinantal representations of $C$ is obtained in this way.

\section{Algorithms to obtain linear determinantal representations}\label{Algorithms}
In this section, we describe two algorithms to obtain (symmetric) linear determinantal representations of $C$
from the data of line bundles on $C$.
The following algorithm is an extension of \cite[Algorithm 1]{Ish1}.

\begin{Alg}[{First algorithm to obtain $M$}]\mbox{}\label{Alg: Alg1}
	\begin{description} 
		\item[Input:] a homogeneous polynomial $F = F(X,Y,Z) \in k[X, Y, Z]$ of degree $d$
		defining a smooth plane curve $C \subset \P^2$, 
		and a non-effective line bundle $\mathcal{L}$ on $C$ of degree $g-1 = d(d-3)/2$.
		\item[Output:] a linear determinantal representation $M$
		corresponding to $\mathcal{L}$ by Theorem \ref{Bij} (1).
		\begin{description}
    		\item[Step 1 (Global sections)] Compute a $k$-basis $\{v_i\}_{1 \le i \le d}$
    		of $W_0 = H^0(C, \mathcal{L}(1))$.
    		\item[Step 2 (First syzygy)] Compute a $k$-basis $\{e_i\}_{1 \le i \le d}$
    		of the kernel $W_1$ of the multiplication map
    		\begin{align*}
    				m \colon H^0(C, \mathcal{L}(1)) \otimes
    				& H^0(C, \mathcal{O}_C(1))
    				\to
    				H^0(C, \mathcal{L}(2)).
    		\end{align*}
    		\item[Step 3 (Output matrix)] 
    		Write the $k$-basis $\{e_i\}_{1 \le i \le d}$ as
    		\[
    			e_i = \sum_{j=1}^d v_j \otimes l_{i, j}(X, Y, Z),
    		\]
    		where $l_{i,j}(X, Y, Z) \in H^0(C, \mathcal{O}_C(1))$ 
    		are $k$-linear forms.
    		Output the matrix \[M = (l_{i,j}(X, Y, Z))_{1 \le i,j \le d}.\]
    	\end{description}
	\end{description}
\end{Alg}

Note that even if $\mathcal{L}$ is a theta characteristic, the matrix $M$ obtained by Algorithm \ref{Alg: Alg1} is not necessarily symmetric.
To obtain a symmetric linear determinantal representation, we take another approach: 
we compute the adjugate matrix of the linear determinantal representation.
This is a classical approach; refer \cite[Section 2]{PSV11}, \cite[Section 6]{Bea77}.

\begin{Alg}[{Second algorithm to obtain $M$}]\mbox{}\label{Alg: Alg2}
	\begin{description}
		\item[Input:] a homogeneous polynomial $F = F(X,Y,Z) \in k[X, Y, Z]$ of degree $d$
		defining a smooth plane curve $C \subset \P^2$, 
		and a non-effective line bundle $\mathcal{L}$ on $C$ of degree $g-1 = d(d-3)/2$.
		\item[Output:] a linear determinantal representation $M$
		corresponding to $\mathcal{L}$ by Theorem \ref{Bij} (1).
		\begin{description}
    		\item[Step 1 (Global sections)] Compute a $k$-basis $\{v_i\}_{1 \le i \le d}$ of $H^0(C, \mathcal{L}(1))$ and a $k$-basis  $\{w_i\}_{1 \le i \le d}$
    		of $H^0(C, \mathcal{L^{\vee}}(1))$.
    		(Here the line bundle $\mathcal{L}^\vee$ is defined as in Theorem \ref{Bij} (1).)
    		\item[Step 2 (Compute adjugate matrices)] For each $1 \le i,j \le d$, compute the image 
    		\[
    		    m(v_iw_j) \in H^0(C, \mathcal{O}_C(d-1))
    		\]
    		of the multiplication map
    		\begin{align*}
    				m \colon H^0(C, \mathcal{L}(1)) &\otimes
    				 H^0(C, \mathcal{L^{\vee}}(1)) \\
    				 & \quad \to
    				H^0(C, \mathcal{O}_C(d-1)).
    		\end{align*}
    		Obtain a matrix $M_a := (m(v_iw_j))_{1 \le i, j \le d}$.
    		\item[Step 3 (Output Matrix)] 
    		It can be shown that all of the entries of the adjugate matrix $M'_a$ of $M_a$ are polynomials of degree $(d-1)^2$ divisible by $F^{d-2}$.
    		Take the adjugate matrix $M'_a$ of the matrix $M_{a}$.
    		Output the matrix \[M := F^{-(d-2)} M'_a.\]
    	\end{description}
	\end{description}
\end{Alg}

To obtain a symmetric linear determinantal representation corresponding to a non-effective theta characteristic $\mathcal{L}$ by Theorem \ref{Bij} (2),
we choose an isomorphism $\mathcal{L^\vee} \cong \mathcal{L}$, take $w_j = v_j$ for each $1 \le j \le d$, and apply Algorithm \ref{Alg: Alg2}. Then the resulting matrix $M$ is symmetric because
the matrix $M_a = (m(v_iv_j))_{1 \le i,j \le d}$ is symmetric.
This algorithm is justified by a similar argument to \cite[6.23.3]{Bea77}.

\subsection{Remarks on computation}
In the following section, we study the problem of computing representatives of all equivalence classes of linear determinantal representations of a smooth plane curve.
We often have the data of $k$-rational divisor classes on $C$,
instead of the data of line bundles on $C$ (or equivalently, $k$-rational divisors).
Recall the exact sequence
\begin{equation*}\label{Eq: obst}
	\xymatrix@C=10pt{
		0 \ar[r] & \Pic(C) \ar[r] & 
		\mathcal{P}ic_{C/k}(k) \ar[r]^-{\delta} & \Br(k) \ar[r] & \Br(k(C)).
	}
\end{equation*}
Here $\mathcal{P}ic_{C/k}$ is the Picard scheme representing the relative Picard functor.
Its identity component $\mathcal{P}ic^0_{C/k}$ is isomorphic to the Jacobian variety $\mathrm{Jac}(C)$.
Because of this sequence, a $k$-rational divisor class $A  \in \mathcal{P}ic_{C/k}(k)$ comes from a line bundle on $C$ if and only if
the Brauer obstruction $\delta(A)$ vanishes in $\Br(k)$ (see \cite[Section 3]{II2}).

When $C$ has a $k$-rational point,
the obstruction homomorphism $\delta$ is the zero map, and every $k$-rational divisor class
comes from a $k$-rational divisor \cite[Proposition 3.2]{II2}.
Moreover, using Gr\"obner basis, we can compute a $k$-rational divisor 
representing a $k$-rational divisor class,
and check whether the divisor is effective or not.

In the following examples, we can easily find rational points, 
hence we can apply our algorithms to each element of the Picard scheme.
We plan to study linear determinantal representations of smooth plane curves without rational points in the future (see \cite[Example 9.5]{Ish2}, \cite[Example 8]{Ish1} for the case of cubics without rational points).
 
%


\section{Examples}\label{Examples}
In this section, we apply our algorithms to two special quartics, the Klein quartic and the Fermat quartic over $\Q$.
We give representatives of all equivalence classes of
linear determinantal representations of these curves over $\Q$.
We use symmetries of these curves to describe results simply.
Note that the linear equivalence class of a divisor $D$ is denoted by $[D]$.

\subsection{Example 1: the Klein quartic over $\Q$}
We apply our algorithms to the Klein quartic $\mathrm{Kl}$ over $\Q$, which is the plane quartic defined by
\[
	X^3Y + Y^3Z + Z^3X = 0.
\]
There are only three $\Q$-rational points on this curve, 
\[
	P_1 = [1:0:0], \quad P_2 = [0:1:0], \quad P_3 = [0:0:1],
\]
and only two points over quadratic fields (see \cite{Tze04})
\[
	Q_1 = [1 : \zeta_3 : \zeta_3^2], \quad Q_2 = [1 : \zeta_3^2 : \zeta_3].
\]
Here $\zeta_3 = (-1 + \sqrt{-3})/2$ is a primitive third root of unity.
The automorphism group $\Aut_\Q(\mathrm{Kl})$ of $\mathrm{Kl}$ defined over $\Q$ is generated by the cyclic permutation
\begin{align*}
    \theta &\colon [X:Y:Z] \mapsto [Y:Z:X].
\end{align*}

The Mordell--Weil group $\mathrm{Jac}({\mathrm{Kl}})(\Q)$ is a finite cyclic group of order $14$ generated by $[D]$, where $[D]$ is the divisor class of the divisor
\begin{align*}
	D &:= P_2 + P_3 - Q_1 - Q_2
\end{align*}
(see \cite{Tze04}).
Note that
\begin{align*}
	\theta([2P_1]) &= [2P_1 + 8D], &
	\theta([D]) &= [-3D].
\end{align*}

The set of $\Q$-rational divisor classes of degree 2 is
\[
    \mathcal{P}ic^2_{\mathrm{Kl}/\Q}(\Q) = \Jac(\mathrm{Kl})(\Q) + [2P_1].
\]
Among them, we find that there exist exactly seven effective $\Q$-rational divisor classes of degree 2 on the Klein quartic. They are:
\begin{align*}
	[2P_1]&, &
	[P_1 + P_2] &= [6D + 2P_1], \\
	[2P_2] &= [12D + 2P_1], &
	[P_1 + P_3] &= [4D + 2P_1], \\
	[P_2 + P_3] &= [10D + 2P_1], &
	[2P_3] &= [8D + 2P_1], \\
	[Q_1 + Q_2] &= [9D + 2P_1].
\end{align*}
Thus there exist exactly seven non-effective $\Q$-rational divisor classes of degree 2.
Their representatives under $\Aut_\Q(\mathrm{Kl})$-action are
\begin{align*}
	[2D + 2P_1] & \quad (\mbox{theta characteristic; fixed by } \theta), \\
	[D + 2P_1],& \\
	[5D + 2P_1] & \quad (\mbox{corresponds to }\mathcal{O}_C(D + 2P_1)^\vee).
\end{align*}

By Algorithm \ref{Alg: Alg1} and Algorithm \ref{Alg: Alg2}, we see that these divisor classes correspond respectively to
\begin{align*}
    M & :=
    \begin{pmatrix}
    X & 0 & 0 & Y \\
    0 & Y & 0 & Z \\
    0 & 0 & Z & X \\
    Y & Z & X & 0
    \end{pmatrix}, \\
    \quad N_{X,Y,Z} &:=
    \begin{pmatrix}
    X & 0 & Y & Y \\
    -Z & -Y+Z & Y-Z & Y \\
    0 & X & -Z & Y-Z \\
    0 & Y & X-Y & -Z
    \end{pmatrix}, \\
    \transp{N_{X,Y,Z}} &\quad (\mbox{the transpose of }N_{X,Y,Z}).
\end{align*}

In conclusion, we have the following theorem.

\begin{thm}
	The Klein quartic $\mathrm{Kl}$ over $\Q$ admits exactly seven equivalence classes of linear determinantal representations over $\Q$.
	They are represented by 
	\begin{align*}
    	M,\quad &
    	N_{X,Y,Z}, \quad N_{Y,Z,X}, \quad N_{Z,X,Y}, \\
    	&\transp{N_{X,Y,Z}}, \quad \transp{N_{Y,Z,X}}, \quad
    	\transp{N_{Z,X,Y}}.
	\end{align*}
	Among them, $M$ gives the only equivalence class of symmetric linear determinantal representations.
\end{thm}

The uniqueness of symmetric linear determinantal representations of the
Klein quartic over $\Q$ was proved by the first and second authors; see \cite[Theorem 1.4]{II1}.

\subsection{Example 2: the Fermat quartic over $\Q$}
Next we apply our algorithms to the Fermat quartic $F_4$ over $\Q$, which is the plane quartic defined by
\[
	X^4 + Y^4 - Z^4 = 0.
\]
See \cite{IIO18, Faddeev} for details of the results in this subsection.

The points on $F_4$ defined over quadratic fields are exhausted by the following sixteen points \cite{Faddeev}:
\begin{align*}
	A_i &= [0:(\sqrt{-1})^i:1] \quad (0 \le i \le 3), \\
	B_j &= [(\sqrt{-1})^j:0:1] \quad (0 \le j \le 3), \\
	E_{1, \pm, \pm} &= \left[ \pm (-1 + \sqrt{-7})/2 : \pm (-1 - \sqrt{-7})/2 : 1 \right],  \\
	E_{2, \pm, \pm} &= \left[ \pm (-1 - \sqrt{-7})/{2} : \pm (-1 + \sqrt{-7})/{2} : 1 \right].
\end{align*}
Among them, there exist exactly four $\Q$-rational points $A_0, A_2, B_0$, and $B_2$.
The points $A_i$ and $B_j$ are hyperflexes, where the tangents meet $F_4$ with multiplicity four.
The automorphism group $\Aut_\Q(F_4)$ is generated by
\begin{align*}
	\theta_1 &\colon [X:Y:Z] \mapsto [-X:Y:Z], \\
	\theta_2 &\colon [X:Y:Z] \mapsto [Y:X:Z].
\end{align*}
They satisfy the relation \(\theta_1^2 = \theta_2^2 = (\theta_2\theta_1)^4 = 1\).

The Mordell--Weil group $\mathrm{Jac}({F_4})(\Q)$ of 
the Jacobian variety of the Fermat quartic $F_4$ is equal to
\begin{align*}\label{Eq: F4}
	(\Z / 4\Z) [D_1] \oplus (\Z / 4\Z) [D_2] \oplus (\Z / 2\Z) [D_3],
\end{align*}
where
\begin{align*}
	D_1 &= A_2 - B_0,\qquad
	D_2 = B_2 - B_0, \\
	D_3 &= A_0 + B_0 - A_2 - B_2. &
\end{align*} 
The $\Aut_\Q(F_4)$-action on $\Jac(F_4)(\Q)$ is described by
\begin{align*}
	\theta_1([D_1]) &= [D_1 - D_2], & \theta_2([D_1]) &= [-D_1-D_3], \\
	\theta_1([D_2]) &= [- D_2], & \theta_2([D_2]) &= [-D_2-D_3], \\
	\theta_1([D_3]) &= [2D_2 + D_3], & \theta_2([D_3]) &= [D_3].
\end{align*}

The set of $\Q$-rational divisor classes of degree 2 is
\[
    \mathcal{P}ic^2_{F_4/\Q}(\Q) = \Jac(F_4)(\Q) + [2B_0].
\]
Among them, we find that there exist exactly sixteen effective $\Q$-rational divisor classes of degree 2 on $F_4$. Ten of them are sums of two $\Q$-rational points on $F_4$:
\begin{align*}
	&[2B_0],&
	&[A_2 + B_0],& 
	&[2A_2],&
	&[B_0 + B_2], \\
	&[A_2 + B_2], &
	&[2B_2],&
	&[A_0 + B_0], &
	&[A_0 + A_2],\\
	&[A_0 + B_2],&
	&[2A_0],&&&&
\end{align*}
and six of them are sums of conjugate pairs of points on $F_4$ defined over quadratic fields:
\begin{align*}
	&[A_1 + A_3],& &[B_1 + B_3],  \\
	&[E_{1, +, +} + E_{2, +, +}],& &[E_{1, +, -} + E_{2, +, -}], \\
	&[E_{1, -, +} + E_{2, -, +}],& &[E_{1, -, -} + E_{2, -, -}].
\end{align*}
Thus there exist exactly sixteen non-effective $\Q$-rational divisor classes of degree 2 on $F_4$.
Their representatives under $\Aut_\Q(F_4)$-action and the stabilizer groups in $\Aut_\Q(F_4)$ are
\begin{align*}
	&[D_3 + 2B_0] \quad &(\mathrm{Stab} &= \{ \id, \theta_1 \}), \\
	&[D_1 - D_2 + 2B_0] \quad &(\mathrm{Stab} &= \{ \id, \theta_2 \}), \\
	&[-D_1 + D_2 + 2B_0] \quad &(\mathrm{Stab} &= \{ \id, \theta_2 \}), \\
	&[2D_1 + D_2 + 2B_0] \quad &(\mathrm{Stab} &= \{ \id, \theta_1 \theta_2\theta_1 \}).
\end{align*}

By Algorithm \ref{Alg: Alg1} and Algorithm \ref{Alg: Alg2}, we see that these divisor classes correspond respectively to
\begin{align*}
    A_{X,Y,Z} :=& 
    \begin{pmatrix}
          X + Z &        -Y &         0 &         Y \\
             -Y &      -X+Z &         0 &         Y \\
              0 &         0 &      -Y+Z &        -X \\
              Y &         Y &        -X &      -Y-Z
    \end{pmatrix}, \\
    B_{X,Y,Z} & :=
    \begin{pmatrix}
     X &  -Z &    Z & Y + Z \\
    -Y &  X + Z &    0 &   0 \\
     0 &  -Y &  X - Z &   0 \\
    -Z &   Z & - Y + Z &   X
    \end{pmatrix}, \\
    \transp{B_{X,Y,Z}} &\quad (\mbox{the transpose of }B_{X,Y,Z}), \\
    C_{X,Y,Z} & :=
    \begin{pmatrix}
     X + Z &         - Z &         0 &  Y + Z \\
         - Y &     X - Z &         0 &         0 \\
             0 & - Y + Z & X + Z &     - Z \\
        - 2Z &             0 &      -Y &  X - Z
    \end{pmatrix}.
\end{align*}

In conclusion, we have the following theorem:

\begin{thm}
	The Fermat quartic $F_4$ over $\Q$ admits exactly sixteen equivalence classes of linear determinantal representations over $\Q$.
	They are represented by 
	\begin{align*}
    	&A_{X, \pm Y, Z}, \quad A_{Y, \pm X, Z} \quad B_{\pm X,\pm Y,Z}, \quad\transp{B_{\pm X,\pm Y,Z}}, \\
		&C_{X,\pm Y,Z}, \quad C_{Y,\pm X,Z}.
	\end{align*}
	Among them, the four equivalence classes represented by $A_{X, \pm Y, Z},\ A_{Y, \pm X, Z}$ give
	the only equivalence classes of symmetric linear determinantal representations.
\end{thm}

%
%

\section*{acknowledgments}
The work of the first author was supported by JSPS KAKENHI Grant Number 13J01450 and 16K17572.
The work of the second author 
was supported by JSPS KAKENHI Grant Number 20674001 and 26800013.
The work of the third author
was supported by JSPS KAKENHI Grant Number 26800011.
Most of calculations 
were done with the aid of  the computer algebra systems 
Maxima \cite{Maxima}, Sage \cite{Sage}, and Singular \cite{Singular, Singular:DivisorsLib}.

\end{document}